\documentclass[a4paper,11pt]{amsart}
\usepackage{graphicx}
\usepackage{mathptmx}
\usepackage{amsmath}
\usepackage{amssymb}
\usepackage{enumitem}
\usepackage{xcolor}
\usepackage{xparse}
\NewDocumentCommand{\eulerian}{omm}
 {%
  \genfrac<>{0pt}{}{#2}{#3}%
  \IfValueT{#1}{_{\!#1}}%
 }

 \newmuskip\pFqmuskip

\newcommand*\pFq[6][8]{%
  \begingroup 
  \pFqmuskip=#1mu\relax
  \mathchardef\normalcomma=\mathcode`,
  \mathcode`\,=\string"8000
  \begingroup\lccode`\~=`\,
  \lowercase{\endgroup\let~}\pFqcomma
  {}_{#2}F_{#3}{\left(\genfrac..{0pt}{}{#4}{#5}\bigg|#6\right)}%
  \endgroup
}
\newcommand{\pFqcomma}{{\normalcomma}\mskip\pFqmuskip}

\newtheorem{theorem}{Theorem}

\begin{document}

\title[A New Approach to Bell and poly-Bell Numbers and Polynomials]{A New Approach to Bell and poly-Bell  Numbers and Polynomials}

\author{DAE SAN KIM}
\address{Department of Mathematics, Sogang University, Seoul 121-742, Republic of Korea}
\email{dskim@sogang.ac.kr}

\author{D. V. Dolgy}
\address{Hanrimwon, Kwangwoon University, Seoul 139-701, Republic of Korea, Institute of Natural Sciences, Far Eastern Federal University, Vladivostok 690950, Russia}
\email{d\_dol@mail.ru }

\author{Hye Kyoung Kim}
\address{Department Of Mathematics Education, Daegu Catholic University, Gyeongsan 38430, Republic of Korea}
\email{hkkim@cu.ac.kr}

\author{Hyunseok  Lee}
\address{Department of Mathematics, Kwangwoon University, Seoul 139-701, Republic of Korea}
\email{luciasconstant@kw.ac.kr}

\author{Taekyun  Kim}
\address{Department of Mathematics, Kwangwoon University, Seoul 139-701, Republic of Korea}
\email{tkkim@kw.ac.kr}


\maketitle

\begin{abstract}
The aim of this paper is to introduce Bell polynomials and numbers of the second kind and poly-Bell polynomials and numbers of the second kind, and to derive their explicit expressions, recurrence relations and some identities involving those polynomials and numbers. We also consider degenerate versions of those polynomials and numbers, namely degenerate Bell polynomials and numbers of the second kind and degenerate poly-Bell polynomials and numbers of the second kind, and deduce their similar results.
\end{abstract}

\section{Introduction}
There are various ways of studying special numbers and polynomials, to mention a few, generating functions, combinatorial methods, $p$-adic analysis, umbral calculus, differential equations, probability theory, special functions and analytic number theory. \par
The aim of this paper is to introduce several special polynomials and numbers, and to study their explicit expressions, recurrence relations and identities involving those polynomials and numbers by using generating functions. \par
Indeed, we introduce Bell polynomials and numbers of the second kind (see \eqref{15}, \eqref{17}) and poly-Bell polynomials and numbers of the second kind (see \eqref{36}). The generating function of Bell numbers of the second kind is the compositional inverse of the generating function of Bell numbers minus the constant term. Then Bell polynomials of the second kind are natural extensions of those numbers. The poly-Bell polynomials of the second kind, which are defined with the help of polylogarithm, become the Bell polynomials of the second kind up to sign when the index of the polylogarithm is $k=1$. \par
We also consider degenerate versions of those numbers and polynomials, namely degenerate Bell numbers and polynomials of the second (see \eqref{28}, \eqref{30}) and degenerate poly-Bell numbers and polynomials (see \eqref{40}), and derive similar results. It is worthwhile to note that degenerate versions of many special numbers and polynomials have been explored in recent years with aforementioned tools and many interesting arithmetical and combinatorial results have been obtained (see [8,9,12,13,18]). In fact, studying degenerate versions can be done not only for polynomials and numbers but also for transcendental functions like gamma functions. For the rest of this section, we recall the necessary facts that are needed throughout this paper. \par                                                                                                                                

The Stirling numbers of the first kind, $S_{1}(n,k)$, are given by
\begin{equation}
\frac{1}{k!}\big(\log(1+t)\big)^{k}=\sum_{n=k}^{\infty}S_{1}(n,k)\frac{t^{n}}{n!},\quad (k\ge 0),\quad (\mathrm{see}\ [5,17]),\label{1}	
\end{equation}
As the inversion formula of \eqref{1}, the Stirling numbers of the second kind, $S_{2}(n,k),$ are given by
\begin{equation}
\frac{1}{k!}\big(e^{t}-1\big)^{k}=\sum_{n=k}^{\infty}S_{2}(n,k)\frac{t^{n}}{n!},\quad (k\ge 0),\quad (\mathrm{see}\ [1,7-14]).\label{2}	
\end{equation}
It is well known that the Bell polynomials are defined as 
\begin{equation}
\mathrm{Bel}_{n}(x)=\sum_{k=0}^{n}S_{2}(n,k)x^{k},\quad (n\ge 0),\quad (\mathrm{see}\ [17]). \label{3}
\end{equation}
From \eqref{3}, we note that 
\begin{equation}
e^{x(e^{t}-1)}=\sum_{n=0}^{\infty}\mathrm{Bel}_{n}(x)\frac{t^{n}}{n!},\quad (\mathrm{see}\ [2,5,6,11,19]).\label{4} 	
\end{equation} 
When $x=1$, $\mathrm{Bel}_{n}=\mathrm{Bel}_{n}(1),\ (n\ge 0)$ are called the Bell numbers. \par 
For any $\lambda\in\mathbb{R}$, the degenerate exponential function is given by 
\begin{equation}
e_{\lambda}^{x}(t)=\sum_{n=0}^{\infty}\frac{(x)_{n,\lambda}}{n!}t^{n}	,\quad (\mathrm{see}\ [3,4,5,18,19]),\label{5}
\end{equation}
where $(x)_{0,\lambda}=1,\ (x)_{n,\lambda}=x(x-\lambda)\cdots(x-(n-1)\lambda),\ (n\ge 1)$. \par
\noindent When $x=1$, we write $e_{\lambda}(t)=e_{\lambda}^{1}(t)$. \par 
The degenerate Stirling numbers of the first kind are defined by 
\begin{equation}
\frac{1}{k!}\big(\log_{\lambda}(1+t)\big)^{k}=\sum_{n=k}^{\infty}S_{1,\lambda}(n,k)\frac{t^{n}}{n!},\ (k\ge 0),\ (\mathrm{see}\ [9]), \label{6}	
\end{equation}
where 
\begin{equation}
\log_{\lambda}(1+t)=\sum_{n=1}^{\infty}\lambda^{n-1}(1)_{n,1/\lambda}\frac{t^{n}}{n!},\quad (\mathrm{see}\ [9]).\label{7}	
\end{equation}
In view of \eqref{2}, the degenerate Stirling numbers of the second kind are defined by 
\begin{equation}
\frac{1}{k!}\big(e_{\lambda}(t)-1\big)^{k}=\sum_{n=k}^{\infty}S_{2}(n,k)\frac{t^{n}}{n!},\quad (\mathrm{see}\ [9]). \label{8}	
\end{equation}
In [11], the degenerate Bell polynomials are defined by 
\begin{equation}
e_{\lambda}^{x}\big(e_{\lambda}(t)-1\big)=\sum_{n=0}^{\infty}\mathrm{Bel}_{n,\lambda}(x)\frac{t^{n}}{n!}.\label{9}
\end{equation}
When $x=1$, $\mathrm{Bel}_{n,\lambda}=\mathrm{Bel}_{n,\lambda}(1),\ (n\ge 0)$, are called the degenerate Bell numbers. \par 
From \eqref{8} and \eqref{9}, we note that 
\begin{equation}
\mathrm{Bel}_{n,\lambda}(x)=\sum_{k=0}^{n}S_{2,\lambda}(n,k)(x)_{k,\lambda},\ (n\ge 0),\ (\mathrm{see}\ [11]).\label{10}	
\end{equation}
The polylogarithm of index $k$ is given by 
\begin{equation}
\mathrm{Li}_{k}(x)=\sum_{n=1}^{\infty}\frac{x^{n}}{n^{k}},\quad (k\in\mathbb{Z},\ |x|<1),\quad(\mathrm{see}\ [1,7,8,10,15]).\label{11}	
\end{equation}
Note that $\mathrm{Li}_{1}(x)=-\log(1-x)$. \par 
Recently, the degenerate polylogarithm is defined as 
\begin{equation}
\mathrm{Li}_{k,\lambda}(x)=\sum_{n=1}^{\infty}\frac{(-\lambda)^{n-1}(1)_{n,1/\lambda}}{(n-1)!n^{k}}x^{n},\quad(|x|<1,\ k\in\mathbb{Z}),\quad(\mathrm{see}\ [9]).\label{12}
\end{equation}
Note that $\mathrm{Li}_{1,\lambda}(x)=-\log_{\lambda}(1-x)$. \par

\section{Bell Polynomials of the second kind} 
From \eqref{4}, we note that 
\begin{displaymath}
e^{x(e^{t}-1)}=\sum_{n=0}^{\infty}\mathrm{Bel}_{n}(x)\frac{t^{n}}{n!} 	
\end{displaymath}
Let $x=1$. Then we have 
\begin{equation}
e^{e^{t}-1}-1	=\sum_{n=1}^{\infty}\mathrm{Bel}_{n}\frac{t^{n}}{n!}. \label{13}
\end{equation}
Let $f(t)=e^{e^{t}-1}-1$. Then the compositional inverse of $f(t)$ is given by 
\begin{equation}
f^{-1}(t)=\log\big(1+\log(1+t)\big).\label{14}
\end{equation}
We consider the new type Bell numbers, called {\it{Bell numbers of the second kind}}, defined by
\begin{equation}
\log\big(1+\log(1+t)\big)=\sum_{n=1}^{\infty}\mathrm{bel}_{n}\frac{t^{n}}{n!}. \label{15}
\end{equation}
Now, we observe that
\begin{align}
&\log\big(1+\log(1+t)\big)=\sum_{k=1}^{\infty}\frac{(-1)^{k-1}}{k}\big(\log(1+t)\big)^{k}\label{16}\\
&=\sum_{k=1}^{\infty}(-1)^{k-1}(k-1)!\frac{1}{k!}\big(\log(1+t)\big)^{k}	=\sum_{k=1}^{\infty}(-1)^{k-1}(k-1)!\sum_{n=k}^{\infty}S_{1}(n,k)\frac{t^{n}}{n!} \nonumber \\
&=\sum_{n=1}^{\infty}\bigg(\sum_{k=1}^{n}(-1)^{k-1}(k-1)!S_{1}(n,k)
\bigg)\frac{t^{n}}{n!}.\nonumber
\end{align}
Therefore, by \eqref{15} and \eqref{16}, we obtain the following theorem. 
\begin{theorem}
For $n\ge 1$, we have 
\begin{displaymath}
(-1)^{n-1}\mathrm{bel}_{n}=\sum_{k=1}^{n}(k-1)!{n \brack k},	
\end{displaymath}
where ${n \brack k} $ are the unsigned Stirling numbers of the first kind. 
\end{theorem}
Also, we consider the new type Bell polynomials, called {\it{Bell polynomials of the second kind}}, defined by
\begin{equation}
\mathrm{bel}_{n}(x)=\sum_{k=1}^{n}(-1)^{k-1}(k-1)!S_{1}(n,k)x^{k},\quad (n\ge 1). \label{17}
\end{equation}
From \eqref{17}, we can derive the following equation. 
\begin{align}
\sum_{n=1}^{\infty} \mathrm{bel}_{n}(x)\frac{t^{n}}{n!}&=\sum_{n=1}^{\infty}\bigg(\sum_{k=1}^{n}(-1)^{k-1}(k-1)!S_{1}(n,k)x^{k}\bigg)\frac{t^{n}}{n!} \label{18} \\
&=\sum_{k=1}^{\infty}(-1)^{k-1}(k-1)!x^{k}\sum_{n=k}^{\infty}S_{1}(n,k)\frac{t^{n}}{n!} \nonumber \\
&=\sum_{k=1}^{\infty}\frac{(-1)^{k-1}k!}{k}x^{k}\frac{1}{k!}\big(\log(1+t)\big)^{k} \nonumber\\
&=\sum_{k=1}^{\infty}\frac{(-1)^{k-1}}{k}x^{k}\big(\log(1+t)\big)^{k}=\log (1+x\log(1+t)).\nonumber
\end{align}
Thus the generating function of Bell polynomials of the second kind is given by 
\begin{equation}
\log\big(1+x\log(1+t)\big)=\sum_{n=1}^{\infty}\mathrm{bel}_{n}(x)\frac{t^{n}}{n!}. \label{19}
\end{equation}
Note here that $\mathrm{bel}_{n}=\mathrm{bel}_{n}(1)$.
From \eqref{19}, we note that 
\begin{equation}
\frac{x}{\big(1+x\log(1+t)\big)(1+t)}=\frac{d}{dt}\log\big(1+x\log(1+t)\big)=\sum_{n=0}^{\infty}\mathrm{bel}_{n+1}(x)\frac{t^{n}}{n!} \label{20}.
\end{equation}
Replacing $t$ by $e^{t}-1$ in \eqref{20}, we get 
\begin{align}
\frac{x}{1+xt}e^{-t}\ &=\ \sum_{k=0}^{\infty}\mathrm{bel}_{k+1}(x)\frac{1}{k!}\big(e^{t}-1\big)^{k} \label{21} \\
&=\ \sum_{k=0}^{\infty}\mathrm{bel}_{k+1}(x)\sum_{n=k}	^{\infty}S_{2}(n,k)\frac{t^{n}}{n!} \nonumber \\
&=\ \sum_{n=0}^{\infty}\bigg(\sum_{k=0}^{n}\mathrm{bel}_{k+1}(x)S_{2}(n,k)\bigg)\frac{t^{n}}{n!}. \nonumber
\end{align}
Taking $x=-1$ in \eqref{21}, we have 
\begin{equation}
\sum_{n=0}^{\infty}\bigg(\sum_{k=0}^{n}\mathrm{bel}_{k+1}(-1)S_{2}(n,k)\bigg)\frac{t^{n}}{n!}=-\frac{1}{1-t}e^{-t}=-\sum_{n=0}^{\infty}d_{n}\frac{t^{n}}{n!}, \label{22}
\end{equation}
where $d_{n}$ is the derangement number ([13]). \par 
Therefore, by comparing the coefficients on both sides of \eqref{22}, we obtain the following theorem. 
\begin{theorem}
For $n\ge 0$, we have 
\begin{displaymath}
\sum_{k=0}^{n}\mathrm{bel}_{k+1}(-1)S_{2}(n,k)=-d_{n}.
\end{displaymath}	
\end{theorem}
Replacing $t$ by $e^{e^{t}-1}-1$ in \eqref{15}, we get 
\begin{align}
t&=\sum_{k=1}^{\infty}\mathrm{bel}_{k}\frac{1}{k!}\Big(e^{e^{t}-1}-1\Big)^{k}=\sum_{k=1}^{\infty}\mathrm{bel}_{k}\sum_{j=k}^{\infty}S_{2}(j,k)\frac{1}{j!}	\big(e^{t}-1\big)^{j}\label{23}\\
&=\sum_{j=1}^{\infty}\sum_{k=1}^{j}\mathrm{bel}_{k}S_{2}(j,k)\sum_{n=j}^{\infty}S_{2}(n,k)\frac{t^{n}}{n!} \nonumber \\
&=\sum_{n=1}^{\infty}\bigg(\sum_{j=1}^{n}\sum_{k=1}^{j}\mathrm{bel}_{k}S_{2}(j,k)S_{2}(n,j)\bigg)\frac{t^{n}}{n!}. \nonumber
\end{align}
Thus we obtain following theorem.
\begin{theorem}
For $n\ge 2$, we have 	
\begin{displaymath}
\sum_{j=1}^{n}\sum_{k=1}^{j}\mathrm{bel}_{k}S_{2}(j,k)S_{2}(n,j)=0, \,\,\, \mathrm{and}\,\,\, \mathrm{bel}_{1}=1. 
\end{displaymath}
\end{theorem}
Replacing $t$ by $e^{t}-1$ in \eqref{19}, we get 
\begin{align}
\log(1+xt)\ &=\ \sum_{k=1}^{\infty}\mathrm{bel}_{k}(x)\frac{1}{k!}\big(e^{t}-1\big)^{k} \label{24} \\
&=\ \sum_{k=1}^{\infty}\mathrm{bel}_{k}(x)\sum_{n=k}^{\infty}S_{2}(n,k)\frac{t^{n}}{n!}\nonumber \\
&=\ \sum_{k=1}^{\infty}\bigg(\sum_{k=1}^{n}\mathrm{bel}_{k}(x)S_{2}(n,k)\bigg)\frac{t^{n}}{n!}.\nonumber 
\end{align}
On the other hand, 
\begin{equation}
\log(1+xt)=\sum_{n=1}^{\infty}\frac{(-1)^{n-1}}{n}x^{n}t^{n}. \label{25}	
\end{equation}
Therefore, by \eqref{24} and \eqref{25}, we obtain the following theorem. 
\begin{theorem}
	For $n\ge 1$, we have 
	\begin{displaymath}
	x^{n}=\frac{(-1)^{n-1}}{(n-1)!}\sum_{k=1}^{n}\mathrm{bel}_{k}(x)S_{2}(n,k). 	
	\end{displaymath}
In particular, 
\begin{displaymath}
	1=\frac{(-1)^{n-1}}{(n-1)!}\sum_{k=1}^{n}\mathrm{bel}_{k}S_{2}(n,k).
\end{displaymath}
\end{theorem}

\section{Degenerate Bell Polynomials of the second kind}
From \eqref{3}, we note that 
\begin{equation}
e_{\lambda}\big(e_{\lambda}(t)-1\big)-1=\sum_{n=1}^{\infty}\mathrm{Bel}_{n,\lambda}\frac{t^{n}}{n!}.\label{26}
\end{equation}
Let $f_{\lambda}(t)=e_{\lambda}\big(e_{\lambda}(t)-1\big)-1$. Then the compositional inverse of $f_{\lambda}(t)$ is given by 
\begin{equation}
f_{\lambda}^{-1}(t)=\log_{\lambda}\big(1+\log_{\lambda}(1+t)\big).\label{27}
\end{equation} \par
We consider the new type degenerate Bell numbers, called {\it{degenerate Bell numbers of the second kind}}, defined by
\begin{equation}
\log_{\lambda}\big(1+\log_{\lambda}(1+t)\big)=\sum_{n=1}^{\infty}\mathrm{bel}_{n,\lambda}\frac{t^{n}}{n!}. \label{28}	
\end{equation}
Now, we observe that 
\begin{align}
\log_{\lambda}\big(1+\log_{\lambda}(1+t)\big)\ &=\ \sum_{k=1}^{\infty}\lambda^{k-1}(1)_{k,1/\lambda}\frac{1}{k!}\big(\log_{\lambda}(1+t)\big)^{k} \label{29}\\
&=\ \sum_{k=1}^{\infty}\lambda^{k-1}(1)_{k,1/\lambda}\sum_{n=k}^{\infty}S_{1,\lambda}(n,k)\frac{t^{n}}{n!}. \nonumber \\
&=\ \sum_{n=1}^{\infty}\bigg(\sum_{k=1}^{n}\lambda^{k-1}(1)_{k,1/\lambda}S_{1,\lambda}(n,k)\bigg)\frac{t^{n}}{n!}. \nonumber
\end{align}
Therefore, by \eqref{28} and \eqref{29}, we obtain the following theorem. 
\begin{theorem}
For $n\ge 1$, we have 
\begin{displaymath}
\mathrm{bel}_{n,\lambda}=\sum_{k=1}^{n}\lambda^{k-1}(1)_{k,1/\lambda}S_{1,\lambda}(n,k).
\end{displaymath}
\end{theorem}
Also, we define the {\it{degenerate Bell polynomials of second kind}} by
\begin{equation}
\mathrm{bel}_{n,\lambda}(x)=\sum_{k=1}^{n}\lambda^{k-1}(1)_{k,1/\lambda}S_{1,\lambda}(n,k)x^{k}.\label{30}
\end{equation}
Note that $\mathrm{bel}_{n,\lambda}=\mathrm{bel}_{n,\lambda}(1)$. \par
From \eqref{30}, we note that 
\begin{align}
\sum_{n=1}^{\infty}\mathrm{bel}_{n,\lambda}(x)\frac{t^{n}}{n!}\ &=\ \sum_{n=1}^{\infty}\bigg(\sum_{k=1}^{n}\lambda^{k-1}(1)_{k,1/\lambda}S_{1,\lambda}(n,k)x^{k}\bigg)\frac{t^{n}}{n!} \label{31} \\
&=\ \sum_{k=1}^{\infty}\lambda^{k-1}(1)_{k,1/\lambda}x^{k}\sum_{n=k}^{\infty}S_{1,\lambda}(n,k)\frac{t^{n}}{n!} \nonumber \\
&=\ \sum_{k=1}^{\infty}\lambda^{k-1}(1)_{k,1/\lambda}x^{k}\frac{1}{k!}\big(\log_{\lambda}(1+t)\big)^{k} \nonumber \\
&=\ \log_{\lambda}\big(1+x\log_{\lambda}(1+t)\big). \nonumber	
\end{align}
Thus the generating function of $\mathrm{bel}_{n,\lambda}(x)$ is given by
\begin{equation}
\log_{\lambda}\big(1+x\log_{\lambda}(1+t)\big)=\sum_{n=1}^{\infty}\mathrm{bel}_{n,\lambda}(x)\frac{t^{n}}{n!}.\label{32}	
\end{equation}
Replacing $t$ by $e_{\lambda}(t)-1$ in \eqref{32}, we get 
\begin{align}
\log_{\lambda}(1+xt)\ &=\ \sum_{k=1}^{\infty}\mathrm{bel}_{k,\lambda}(x)\frac{1}{k!}\big(e_{\lambda}(t)-1\big)^{k}\label{33} \\
&=\ \sum_{k=1}^{\infty}\mathrm{bel}_{k,\lambda}(x)\sum_{n=k}^{\infty}S_{2,\lambda}(n,k)\frac{t^{n}}{n!}\nonumber \\
&=\ \sum_{n=1}^{\infty}\bigg(\sum_{k=1}^{n}\mathrm{bel}_{k,\lambda}(x)S_{2,\lambda}(n,k)\bigg)\frac{t^{n}}{n!}. \nonumber
\end{align}
On the other hand, 
\begin{equation}
\log_{\lambda}(1+xt)=\sum_{n=1}^{\infty}\lambda^{n-1}(1)_{n,1/\lambda}x^{n}\frac{t^{n}}{n!}.\label{34}
\end{equation}
Therefore, by \eqref{33} and \eqref{34}, we obtain the following theorem.
\begin{theorem}
For $n\ge 1$, we have 
\begin{displaymath}
x^{n}=\frac{\lambda^{1-n}}{(1)_{n,1/\lambda}}=\sum_{k=1}^{n}\mathrm{bel}_{k,\lambda}(x)S_{2,\lambda}(n,k).
\end{displaymath}
In particular, 
\begin{displaymath}
\lambda^{n-1}(1)_{n,1/\lambda}=\sum_{k=1}^{n}\mathrm{bel}_{k,\lambda}S_{2,\lambda}(n,k).
\end{displaymath}
\end{theorem}
Replacing $t$ by $e_{\lambda}\big(e_{\lambda}(t)-1\big)-1$ in \eqref{28}, we have 
\begin{align}
t\ &= \sum_{k=1}^{\infty}\mathrm{bel}_{k,\lambda}\frac{1}{k!}	\big(e_{\lambda}(e_{\lambda}(t)-1)-1\big)^{k}\label{35} \\
&= \sum_{k=1}^{\infty}\mathrm{bel}_{k,\lambda}\sum_{j=k}^{\infty}S_{2,\lambda}(j,k)\frac{1}{j!}\big(e_{\lambda}(t)-1\big)^{j} \nonumber \\
&=\ \sum_{j=1}^{\infty}\bigg(\sum_{k=1}^{j}\mathrm{bel}_{k,\lambda}S_{2,\lambda}(j,k)\bigg)\sum_{n=j}^{\infty}S_{2,\lambda}(n,j)\frac{t^{n}}{n!} \nonumber \\
&=\ \sum_{n=1}^{\infty}\bigg(\sum_{j=1}^{n}\sum_{k=1}^{j}\mathrm{bel}_{k,\lambda}S_{2,\lambda}(j,k)S_{2,\lambda}(n,j)\bigg)\frac{t^{n}}{n!}. \nonumber 
\end{align}
Therefore, by comparing the coefficients on both sides of \eqref{35}, we obtain the following theorem. 
\begin{theorem}
For $n\ge 2$, we have 	
\begin{displaymath}
\sum_{j=1}^{n}\sum_{k=1}^{j}\mathrm{bel}_{k,\lambda}S_{2,\lambda}(j,k)S_{2,\lambda}(n,j)=0,\,\,\,\mathrm{and}\,\,\, \mathrm{bel}_{1,\lambda}=1. 
\end{displaymath}
\end{theorem}

\section{Poly-Bell Polynomials of the second kind}
Now, we consider the {\it{poly-Bell polynomials of the second kind}} which are defined as
\begin{equation}
\mathrm{Li}_{k}\big(-x\log(1-t)\big)=\sum_{n=1}^{\infty}\mathrm{bel}_{n}^{(k)}(x)\frac{t^{n}}{n!}.\label{36}
\end{equation}
When $x=1$, $\mathrm{bel}_{n}^{(k)}=\mathrm{bel}_{n}^{(k)}(1)$ are called the {\it{poly-Bell numbers of the second kind}}. \par 
From \eqref{11}, we note that 
\begin{align}
\mathrm{Li}_{k}\big(-x\log(1-t)\big)\ &=\ \sum_{l=1}^{\infty}\frac{(-1)^{l}}{l^{k}}x^{l}l!\frac{1}{l!}\big(\log(1-t)\big)^{l}\label{37} \\
&=\ \sum_{l=1}^{\infty}\frac{(-1)^{l}}{l^{k-1}}(l-1)!x^{l}\sum_{n=l}^{\infty}(-1)^{n}S_{1}(n,l)\frac{t^{n}}{n!} \nonumber \\
&=\ \sum_{n=1}^{\infty}\bigg(\sum_{l=1}^{n}\frac{(-1)^{n-l}}{l^{k-1}}(l-1)!x^{l}S_{1}(n,l)\bigg)\frac{t^{n}}{n!}.\nonumber
\end{align}
Therefore, by \eqref{36} and \eqref{37}, we obtain the following theorem. 
\begin{theorem}
For $n\ge 1$, we have 
\begin{displaymath}
\mathrm{bel}_{n}^{(k)}(x)=\sum_{l=1}^{n}\frac{x^{l}}{l^{k-1}}(l-1)!{n \brack l}.	
\end{displaymath}
In particular, 
\begin{displaymath}
\mathrm{bel}_{n}^{(k)}=\sum_{l=1}^{n}\frac{1}{l^{k-1}}(l-1)!{n \brack l}.	
\end{displaymath}
\end{theorem}
Note that 
\begin{displaymath}
\mathrm{bel}_{n}^{(1)}(x)=\sum_{l=1}^{n}x^{l}(l-1)!	{n \brack l}=(-1)^{n-1}\mathrm{bel}_{n}(x). 
\end{displaymath}
Indeed, 
\begin{align*}
\mathrm{Li}_{1}\big(-x\log(1-t)\big)\ &=\ -\log\big(1+x\log(1-t)\big) \\
&=\ \sum_{n=1}^{\infty}\mathrm{bel}_{n}(x)(-1)^{n-1}\frac{t^{n}}{n!}.	
\end{align*}
Replacing $t$ by $1-e^{-t}$ in \eqref{36}, we get 
\begin{align}
\mathrm{Li}_{k}(xt)\ &=\ \sum_{l=1}^{\infty}\mathrm{bel}_{l}^{(k)}(x)\frac{1}{l!}\big(1-e^{-t}\big)^{l} \label{38} \\
&=\ \sum_{l=1}^{\infty}\mathrm{bel}_{l}^{(k)}(x)(-1)^{l}\sum_{n=l}^{\infty}S_{2}(n,l)(-1)^{n}\frac{t^{n}}{n!}.\nonumber \\
&=\ \sum_{n=1}^{\infty}\bigg(\sum_{l=1}^{n}(-1)^{n-l}\mathrm{bel}_{l}^{(k)}(x)S_{2}(n,l)\bigg)\frac{t^{n}}{n!}. \nonumber 	
\end{align}
From \eqref{11} and \eqref{38}, we note that 
\begin{equation}
\frac{x^{n}}{n^{k}}=\frac{1}{n!}\sum_{l=1}^{n}(-1)^{n-l}\mathrm{bel}_{l}^{(k)}(x)S_{2}(n,l).\label{39}
\end{equation}
Therefore, by \eqref{39}, we obtain the following theorem. 
\begin{theorem}
For $n\ge 1$, we have 
\begin{displaymath}
x^{n}=\frac{n^{k-1}}{(n-1)!}\sum_{l=1}^{n}(-1)^{n-l}\mathrm{bel}_{l}^{(k)}(x)S_{2}(n,l).
\end{displaymath}	
\end{theorem}

\section{Degenerate Poly-Bell Polynomials of the Second Kind}
We define the {\it{degenerate poly-Bell polynomials of the second kind}} by
\begin{equation}
\mathrm{Li}_{k,\lambda}\big(-x\log_{\lambda}(1-t)\big)=\sum_{n=1}^{\infty}\mathrm{bel}_{n,\lambda}^{(k)}(x)\frac{t^{n}}{n!}.\label{40}	
\end{equation}
When $x=1$, $\mathrm{bel}_{n,\lambda}^{(k)}=\mathrm{bel}_{n,\lambda}^{(k)}(1)$ are called the {\it{degenerate poly-Bell numbers of the second}}. \par
From \eqref{13}, we note that 
\begin{align}
\mathrm{Li}_{k,\lambda}\big(-x\log_{\lambda}(1-t)\ &=\ \sum_{l=1}^{\infty}\frac{(-\lambda)^{l-1}(1)_{l,1/\lambda}}{(l-1)!l^{k}}\big(-x\log_{\lambda}(1-t)\big)^{l} \label{41} \\
&=\ -\sum_{l=1}^{\infty}\frac{(1)_{l,1/\lambda}}{l^{k-1}}\lambda^{l-1}x^{l}\frac{1}{l!}\big(\log_{\lambda}(1-t)\big)^{l}. \nonumber \\
&= -\sum_{l=1}^{\infty}\frac{(1)_{l,1/\lambda}}{l^{k-1}}\lambda^{l-1}x^{l}\sum_{n=l}^{\infty}S_{1,\lambda}(n,l)(-1)^{n}\frac{t^{n}}{n!}\nonumber \\
&=\ \sum_{n=1}^{\infty}\bigg((-1)^{n-1}\sum_{l=1}^{n}\frac{1}{l^{k-1}}(1)_{l,1/\lambda}\lambda^{l-1}x^{l}S_{1,\lambda}(n,l)\bigg)\frac{t^{n}}{n!}. \nonumber
\end{align}
Therefore, by \eqref{40} and \eqref{41}, we obtain the following theorem. 
\begin{theorem}
For $n\ge 1$, we have 
\begin{displaymath}
(-1)^{n-1}\mathrm{bel}_{n,\lambda}^{(k)}(x)=\sum_{l=1}^{n}\frac{1}{l^{k-1}}(1)_{l,1/\lambda}\lambda^{l-1}x^{l}S_{1,\lambda}(n,l). 	
\end{displaymath}
\end{theorem} 
For $k=1$, we have 
\begin{displaymath}
(-1)^{n-1}\mathrm{bel}_{n,\lambda}^{(1)}(x)=\sum_{l=1}^{n}(1)_{l,1/\lambda}\lambda^{l-1}x^{l}S_{1,\lambda}(n,l)=\mathrm{bel}_{n,\lambda}(x),\quad (n\ge 0).
\end{displaymath}
Indeed, 
\begin{displaymath}
\mathrm{Li}_{1,\lambda}\big(-x\log_{\lambda}(1-t)\big)=-\log_{\lambda}\big(1+x\log_{\lambda}(1-t)\big)=\sum_{n=1}^{\infty}(-1)^{n-1}\mathrm{bel}_{n,\lambda}(x)\frac{t^{n}	}{n!}. 
\end{displaymath}
Replacing $t$ by $1-e_{\lambda}(-t)$ in \eqref{40}, we get 
\begin{align}
\mathrm{Li}_{k,\lambda}(xt)\ &=\ \sum_{l=1}^{\infty}	\mathrm{bel}_{l,\lambda}^{(k)}(x)\frac{1}{l!}\big(1-e_{\lambda}(-t)\big)^{l} \label{42} \\
&=\ \sum_{l=1}^{\infty} \mathrm{bel}_{l,\lambda}^{(k)}(x)(-1)^{l}\frac{1}{l!}\big(e_{\lambda}(-t)-1\big)^{l}\nonumber \\
&=\ \sum_{l=1}^{\infty} \mathrm{bel}_{l,\lambda}^{(k)}(x)(-1)^{l}\sum_{n=l}^{\infty}S_{2,\lambda}(n,l)(-1)^{n}\frac{t^{n}}{n!} \nonumber \\
&=\ \sum_{n=1}^{\infty}\bigg(\sum_{l=1}^{n}(-1)^{n-l} \mathrm{bel}_{l,\lambda}^{(k)}(x)S_{2,\lambda}(n,l)\bigg)\frac{t^{n}}{n!}. \nonumber 
\end{align}
On the other hand, 
\begin{equation}
\mathrm{Li}_{k,\lambda}(xt)=\sum_{n=1}^{\infty}\frac{(-\lambda)^{n-1}(1)_{n,1/\lambda}}{(n-1)!n^{k}}x^{n}t^{n}=\sum_{n=1}^{\infty}\frac{(-\lambda)^{n-1}(1)_{n,1/\lambda}}{n^{k-1}}x^{n}\frac{t^{n}}{n!}.\label{43}
\end{equation}
From \eqref{42} and \eqref{43}, we get the following result.
\begin{theorem}
For $n \ge 1$, we have
\begin{displaymath}
\frac{(-\lambda)^{n-1}(1)_{n,1/\lambda}}{n^{k-1}}x^{n}=\sum_{l=1}^{n}(-1)^{n-l}\mathrm{bel}_{l,\lambda}^{(k)}(x)S_{2,\lambda}(n,l).
\end{displaymath}
\end{theorem}

\section{Conclusion} 

By means of various different tools, degenerate versions of many special polynomials and numbers have been studied in recent years. Here we introduced Bell polynomials of the second kind, poly-Bell polynomials of the second kind and their degenerate versions, namely degenerate Bell polynomials of the second kind and degenerate poly-Bell polynomials of the second kind. By using generating functions, we explored their explicit expressions, recurrence relations and some identities involving those polynomials and numbers.\par
It is one of our future projects to continue this line of research, namely to explore many special numbers and polynomials and their degenerate versions with the help of various different tools.
\vspace{0.2in}

\noindent{\bf{Acknowledgments:}} 

\vspace{0.2in}

\noindent{\bf{Funding:}} 

\vspace{0.2in}

\noindent{\bf {Availability of data and materials:}}
Not applicable.

\vspace{0.1in}

\noindent{\bf {Competing interests:}}
The authors declare no conflict of interest.

\vspace{0.1in}

\noindent{\bf{Authors' contributions:}} T.K. and D.S.K. conceived of the framework and structured the whole paper; T. K. and D.S.K. wrote the paper; H.K.K., D.V.D. and H.L. completed the revision of the article. All authors have read and agreed to the published version of the manuscript.

\vspace{0.1in}

\noindent{\bf{Author details:}}

\end{document}